\documentclass[final]{ifacconf}

\usepackage[utf8x]{inputenc}
\usepackage{amssymb,amsmath}
\usepackage{mathtools}
\usepackage{mathrsfs}
\usepackage{amssymb}
\usepackage{natbib}           
\usepackage{graphicx}          
\usepackage{url}	    

\usepackage{algorithm}
\usepackage{algorithmicx}
\usepackage{algpseudocode}
\usepackage{float}
\usepackage{algpseudocode}
\usepackage{varwidth}

\usepackage{color,xcolor}

\usepackage{units}

\newcommand{\Laplace}{\mathcal{L}}
\newcommand{\N}{\mathbb{N}}

\renewcommand{\Re}{\mathbb{R}}

\newcommand{\D}{\operatorname{D}}
\newcommand{\rlD}{{_{\mathrm{rl}}\D}}
\newcommand{\glD}{{_{\mathrm{gl}}\D}}
\newcommand{\glDD}{{_{\mathrm{gl}}\Delta}}
\newcommand{\cD}{{_{\mathrm{c}}\D}}
\newcommand{\rlI}{{_{\mathrm{rl}}\operatorname{I}}}

\usepackage{titlesec}
\begin{document}

\begin{frontmatter}

\title{Modeling and administration scheduling of fractional-order 
pharmacokinetic  systems%
\thanksref{footnoteinfo}} 

%
%
\thanks[footnoteinfo]{The work of the last
author was supported by the KU Leuven Research Council under BOF/STG-15-043.}

\author[IMTLUCCA]{Domagoj Herceg},
\author[NTUA]{Sotiris Ntouskas},
\author[KUL]{Pantelis Sopasakis},
\author[PHARMA]{Aris Dokoumetzidis},
\author[PHARMA]{Panos Macheras},
\author[NTUA]{Haralambos Sarimveis}  and
\author[KUL]{Panagiotis Patrinos}

\address[IMTLUCCA]{IMT School for Advanced Studies Lucca, Piazza San Francesco, 19, 55100 Lucca, Italy.}                                              
\address[KUL]{KU Leuven, Department of Electrical Engineering (ESAT), STADIUS Center for Dynamical Systems, Signal Processing and Data Analytics \& Optimization in Engineering (OPTEC), Kasteelpark Arenberg 10, 3001 Leuven, Belgium.}
\address[NTUA]{National Technical University of Athens (NTUA), School of Chemical Engineering, 9 Heroon Polytechneiou Street, 15780 Zografou~Campus, Athens, Greece.}
\address[PHARMA]{Laboratory of Biopharmaceutics \& Pharmacokinetics, Faculty of Pharmacy, National \& Kapodistrian University of Athens, 
Panepistimiopolis, 15771 Zografou, Greece}


\begin{abstract}  
Fractional-order dynamical systems were recently introduced in the 
field of pharmacokinetics where they proved powerful tools for 
modeling the absorption, disposition, distribution and excretion of drugs
which are liable to anomalous diffusion, deep tissue trapping and 
other nonlinear phenomena. In this paper we present several 
ways to simulate such fractional-order pharmacokinetic models and 
we evaluate their accuracy and complexity on a fractional-order 
pharmacokinetic model of Amiodarone, an anti-arrhythmic drug.
We then propose an optimal administration scheduling scheme 
and evaluate it on a population of patients.
\end{abstract}

\end{frontmatter}

\section{Introduction}
Pharmacokinetic (PK) models are systems of differential 
equations  which simulate the dynamic response of living 
organisms in terms of the concentration of a drug or any 
other substance in different compartments (organs) of the 
body following its administration to the body. PK models 
can assist in designing effective and safe administration 
strategies for individual patients or populations thereof. 
Among the different types of PK modeling, 
fractional PK models have attracted the interest of researchers 
in the field because they can model phenomena like anomalous 
diffusion, deep tissue trapping and diffusion across 
fractal manifolds, which traditional PK models fail to describe~\citep{DokMac08,DokMac11,DokMagMac10}.

However, simulating such dynamics is not as straightforward 
as with integer-order systems.
Analytical solutions are rarely available and, even then,
the evaluation of the solution requires a numerical approximation
method~\citep{Kaczorek2011}. 
The availability of accurate discrete-time approximations
of the trajectories of such systems is important not only 
for simulating but also for the design of 
open-loop or closed-loop (such as with model predictive control) 
administration strategies~\citep{SopNtoSar15}. Here,
we compare several approximation methods (such as ones based on 
the Gr\"{u}nwald-Letnikov operator and the Oustaloup filter) 
using a specific PK model taken from the recent 
literature~\citep{DokMagMac10} and we discuss their accuracy
and fitness for control design.

In this paper we provide a survey of different 
approaches for modeling fractional-order systems which 
arise in pharmacokinetics. Our discussion revolves around 
the case study of Amiodarone, an anti-arrhythmic drug that exhibits 
fractional-order dynamics. We identify three major 
classes of numerical algorithms for simulating fractional-order
systems in the literature: 
(i) using rational transfer functions, (ii) time-domain methods
and (iii) the numerical inverse Laplace approach.
We discuss their merits and limitations and we present a
comparative assessment regarding precision of various methods.
Our goal, however, is to single out a method which is most suitable for
controller design. In the last section we formulate an optimal control
problem for administration scheduling to confirm our findings.

\section{Fractional Pharmocokinetics: Modeling and Simulation}

\subsection{Fractional-order pharmacokinetics}
\begin{figure}
 \centering
 \includegraphics[width=0.33\textwidth]{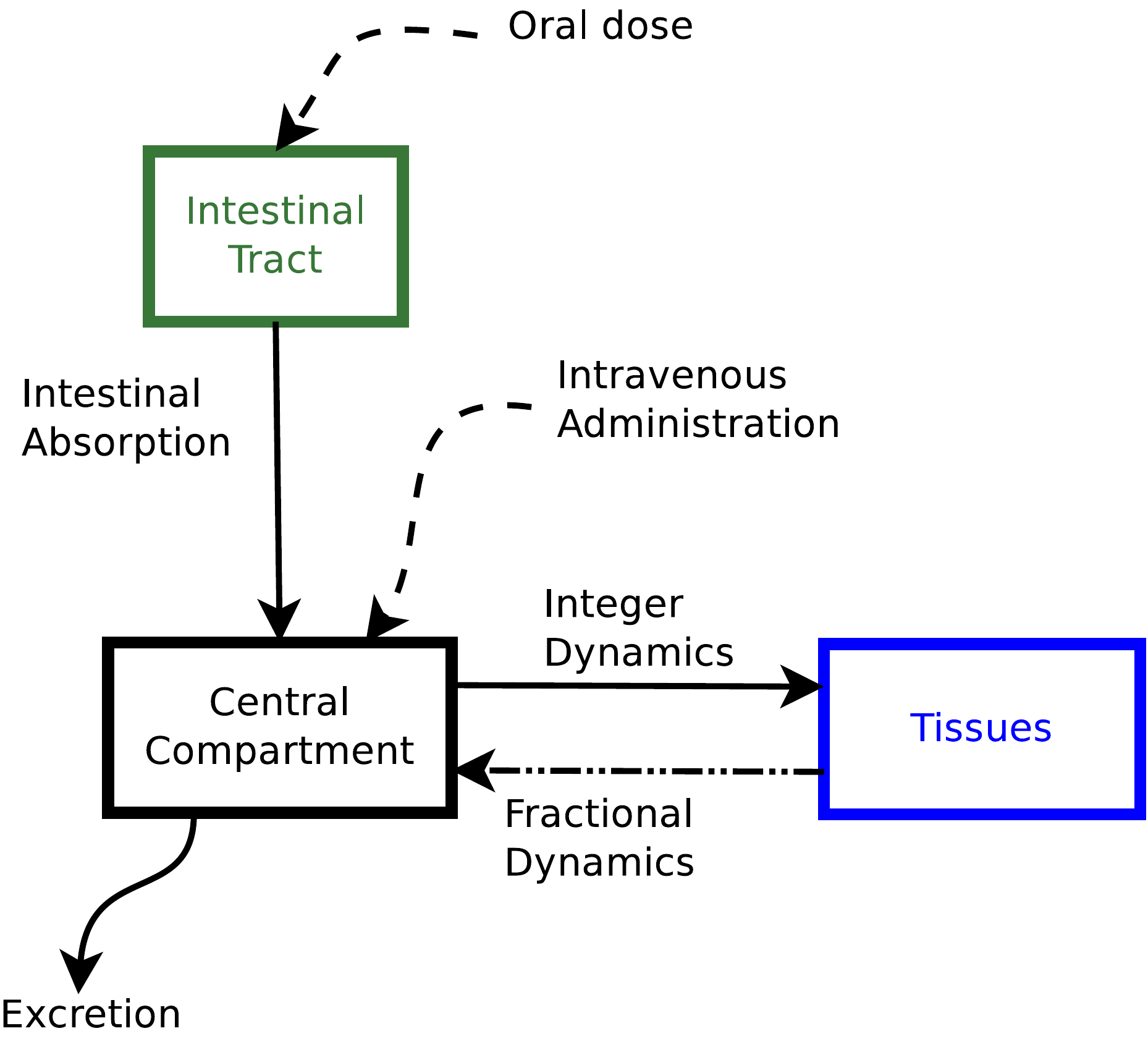}
 \caption{Strucutre of the fractional-order PK model of 
          Amiodarone.}
 \label{fig:pk-topology}          
\end{figure}

Amiodarone is an anti-arrhythmic agent which can be 
administered either intravenously (\emph{i.v.}) 
or orally~\citep{Kuehlkamp1999271}.
It is well-known for its highly nonlinear
non-exponential dynamics and singular long-term
accumulation pattern.
\cite{DokMagMac10} modeled the pharmacokinetic 
distribution of Amiodarone with a fractional compartmental
model following a single intravenous and a single oral
dose. The compartmental topology of the model is presented 
in Figure~\ref{fig:pk-topology} where it is shown that 
the diffusion from the tissues to the central
compartment is governed by a fractional-order dynamics.

Let $A_1$ and $A_2$ be the amounts of Amiodarone (in~$\mathrm{ng}$)
in the plasma and the tissues respectively and $u$ be the 
administration rate (in $\mathrm{ng}/\mathrm{day}$).
We assume that the drug is administered directly into the
central (plasma) compartment while the control objective
is the concentration of the drug in the tissues attains
a prescribed value (set-point).
The fractional dynamical model we employ reads as follows:
\begin{subequations}
\label{eq:pkModelAmiodarone}
 \begin{eqnarray}
  \frac{\mathrm{d}A_1}{\mathrm{d}t}&=&-(k_{12}+k_{10})A_1 + k_{21}\cdot\, \cD^{1-\alpha} A_2 + u,\\
\frac{\mathrm{d}A_2}{\mathrm{d}t}&=& k_{12}A_1 - k_{21}\cdot \, \cD^{1-\alpha} A_2,
 \end{eqnarray}
\end{subequations}
with $\alpha\in(0,1)$ and $\cD^{1-\alpha}$ is the \textit{Caputo
fractional derivative} which is defined in the following 
section.

\subsection{Fractional-order derivatives}
Several fractional-order derivatives have been proposed in the 
literature the most popular of which are the Riemann-Liouville 
$\rlD^\alpha$, the Caputo $\cD^\alpha$ and the 
Gr\"{u}nwald-Letnikov $\glD^\alpha$ derivatives~\citep{Sam+93}. 
These operators are used to formulate fractional-order 
differential equations, that is functional equations of the form
\begin{align}
 F(x(t), \D^{\alpha_1}x(t), \ldots, \D^{\alpha_p}x(t)) = 0,
\end{align}
where $\D^{\alpha}$ is a generalized derivative of order 
$\alpha\geq 0$.
Typically, the Caputo derivative is used in this context as 
the initial conditions are easier to postulate.

The generalized Riemann-Liouville fractional-order
integral operator of order $\alpha>0$ is given by
\begin{equation}
(\rlI^\alpha f)(t)=
 \tfrac{1}{\Gamma(\alpha)}\int_{0}^t(t-\tau)^{\alpha-1}f(\tau)
   \mathrm d \tau,\ t\geq 0.
   \label{eq:RL_frac_integral}
\end{equation}
For $\alpha\in\Re$ let us denote by $m=\lceil \alpha \rceil$ the
smallest natural number $m$ so that $m\geq \alpha$. The following
operator is known as the \textit{Caputo} derivative of order 
$\alpha$:
\begin{equation}
(\cD^\alpha f)(t) = 
  {\rlI}^{\, m - \alpha}\frac{\mathrm d^{m} f(t)}{\mathrm d t^{m}}.
\end{equation}
The Gr\"{u}nwald-Letnikov fractional-order derivative
is defined as 
\begin{align}
 (\glD^\alpha f)(t) = \lim_{h\to 0} 
   \tfrac{1}{h^\alpha}
   \sum_{i=0}^{\infty}(-1)^i \tbinom{\alpha}{i} f(t-ih), 
\end{align}
where ${\alpha \choose 0}=1$ and $\binom{\alpha}{i}=
\prod_{l=0}^{i-1}\tfrac{\alpha-l}{l+1}$.

\begin{figure}[h]
\centering
\includegraphics[scale = 0.42]{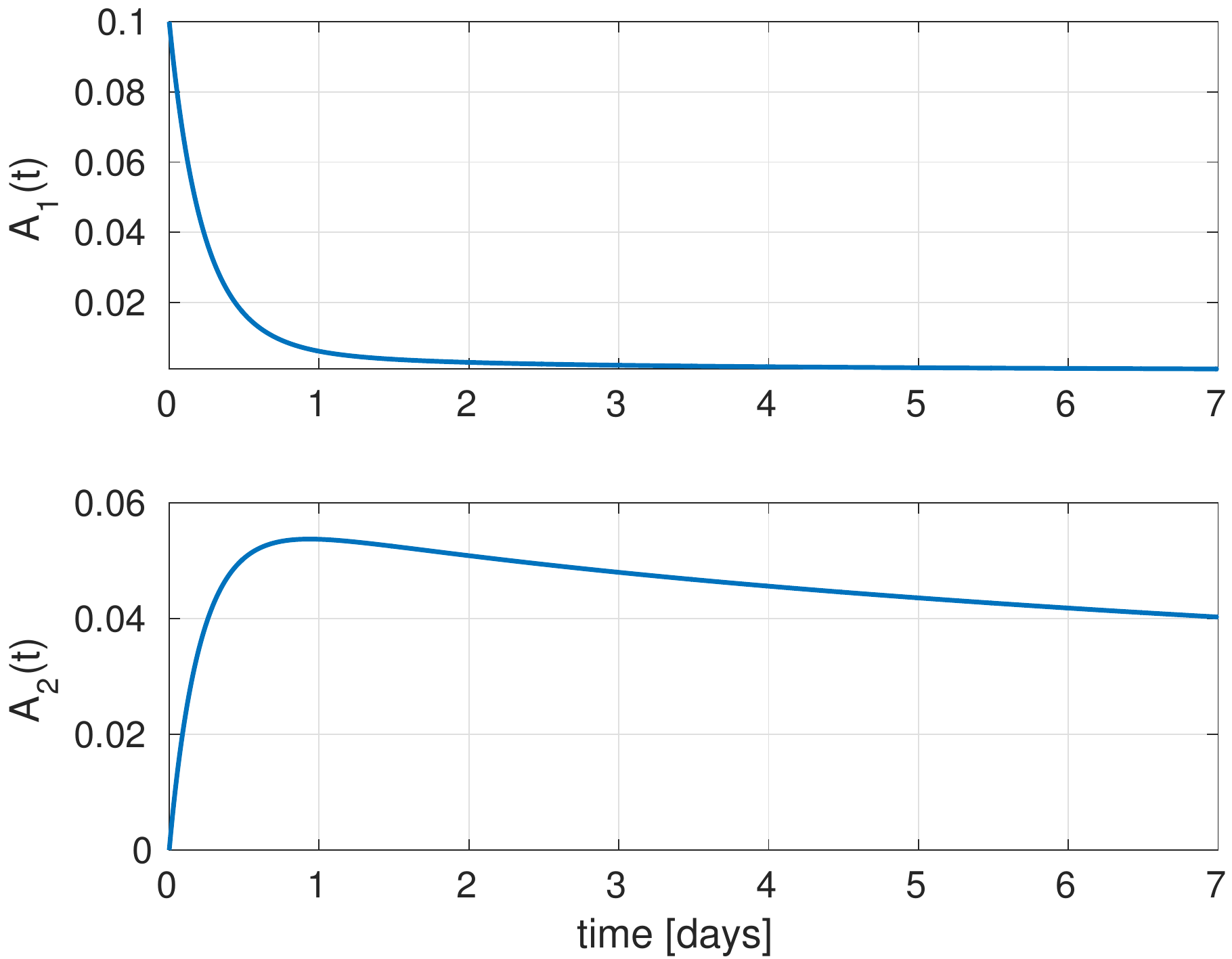}
\caption{Open loop response of system~\eqref{eq:pkModelAmiodarone} with initial conditions $A_1(0) =\unit[0.1]{ng}, A_2(0) =\unit[0]{ng}$. }
\label{fig:open-loop-sim}        
\end{figure}

The Laplace transform of the Caputo derivative of fractional order 
$\alpha\in(0,1)$ with zero initial conditions  is given as 
$
 \Laplace \left[(\cD^\alpha)(t) \right] = s^\alpha F(s)
$, where $F(s)$ is the Laplace transform of function $f(t)$.
The transfer function $G_i(s) = A_i(s)/U(s)$ which
associates the administration rate $U(s)=\Laplace u(t)$
to the concentrations $\hat{A}_i(s)=\Laplace A_i(t)$ are:
\begin{subequations}
\begin{align}
 G_1(s)&=\frac{s^\alpha+k_{21}}
              {s^{\alpha+1}+k_{21}s+(k_{12}+k_{10})s^\alpha+k_{10}k_{21}},\\
 G_2(s)&=\frac{k_{12}s^{\alpha-1}}
              {s^{\alpha+1}+k_{21}s+(k_{12}+k_{10})s^\alpha+k_{10}k_{21}},
\end{align}
\end{subequations}
with $\alpha=0.587$, $k_{10} = \unit[1.4913]{day^{-1}}$,
$k_{12}=\unit[2.9522]{day^{-1}}$ and $k_{21}=\unit[0.4854]{day^{-\alpha}}$.
The concentration of Amiodarone in the two compartments
with initial condition $A_1(0)=\unit[0.1]{ng}$ and $A_2(0)=\unit[0]{ng}$ 
is shown in Figure~\ref{fig:open-loop-sim}.

%
%

\subsection{Solutions of FDEs}
There can be identified four types of solutions for 
fractional-order differential equations: (i) analytical solutions,
(ii) approximations in the $s$-domain using integer-order
rational transfer functions and (iii) numerical approximation 
schemes in the discrete time domain, (iv) the numerical 
inverse Laplace transformation.

\subsubsection{Analytical solutions.}
Analytical solutions, when available,
involve special functions such as the Mittag-Leffler 
function $\mathcal{E}_{\alpha,\beta}(t)=\sum_{k=0}^{\infty}
{t^k}/{\Gamma(\alpha k + \beta)}$ whose evaluation 
requires in turn some numerical approximation scheme. 
Typically for the evaluation of the this function we resort 
to solving an FDE numerically~\citep{Garrappa15}.

\subsubsection{Transfer function approximations.}
Rational approximations aim at approximating the transfer
function of a fractional-order system --- which involves
terms of the form $s^\alpha$ --- by ordinary transfer 
functions of the form 
\begin{align}
 T(s) = \frac{P(s)}{Q(s)},
\end{align}
where $P$ and $Q$ are polynomials and the degree of $P$ is
no larger than the degree of $Q$.

%
%
\textit{Pad\'{e} Approximation:} 
The Pad\'{e} approximation of order $[m/n]$, $m,n\in\N$,
at a point $s_0$ is rather popular and leads to rational functions 
with $\deg P =m$ and $\deg Q = n$~\citep{SILVA2006373}. 

%
%
\textit{Matsuda-Fujii Method:}
This method consists in interpolating a function $H(s)$, which is treated
as a black box, across a set of logarithmically spaced points~\citep{Matsuda1993}. 
By letting the selected points be $s_k$, $k=0,1,2,\dots$, the approximation is 
written as the continued fractions expansion
\begin{equation}
H(s) = \alpha_0 + \frac{s-s_0}{\alpha_1 + \frac{s-s_1}{\alpha_2 + \frac{s-s_2}{\alpha_3+\dots}}}
\end{equation}
where, $\alpha_i = \upsilon_i (s_i)$,  $\upsilon_0(s) = H(s)$, $\upsilon_{i+1}(s) = \frac{s-s_i}{\upsilon_i(s) - \alpha_i}$

%
%
\textit{Oustaloup's method:}
Oustaloup's method is based on the approximation of a function of the form:
\begin{equation}
H(s) = s^{\alpha}, 
\end{equation}
with $\alpha>0$ by a rational function
\begin{equation}\label{eq:oustaloup}
\widehat{H}(s) = c_0 \prod_{k=-N}^N \frac{s+\omega_k}{s+\omega'_k}
\end{equation}
within a range of frequencies from $\omega_b$ to $\omega_h$~\citep{OusLev+00}.
The Oustaloup method offers an approximation at frequencies 
which are geometrically distributed about the characteristic 
frequency $\omega_u = \sqrt{\omega_b\omega_h}$ --- the geometric 
mean of $\omega_b$ and $\omega_h$.
The parameters $\omega_k$ and $\omega_k'$ are determined via the design 
formulas~\citep{Petras11}
\begin{subequations}
 \begin{align}
  \omega_k' = \omega_b \left(\frac{\omega_h}{\omega_b}\right)^{
		\frac{k+N+0.5(1+\alpha)}{2N+1}
	      },\\
  \omega_k = \omega_b \left(\frac{\omega_h}{\omega_b}\right)^{
		\frac{k+N+0.5(1-\alpha)}{2N+1}
	      },\\
  c_0 = \left(\frac{\omega_h}{\omega_b}\right)^{-\frac{r}{2}}
        \prod_{k=-N}^{N}\frac{\omega_k}{\omega_k'}.
 \end{align}
\end{subequations}
Parameters $\omega_b$, $\omega_h$ and $N$ are design parameters 
of the Oustaloup method.

%
%
\textit{Other methods:}
There are a few more methods which have been proposed in the 
literature to approximate fractional-order systems by 
rational transfer functions such as~\citep{Charef,CarHal64}, 
as well as data-driven system identification 
techniques~\citep{Gao2012}.
Nice~\citep{state-space-frac}.

\subsubsection{Time domain approximations.}
Several methods have been proposed which attempt to 
approximate the solution to a fractional-order initial
value problem in the time domain.

%
%
\textit{Gr\"{u}nwald-Letnikov:}
This is the method of choice in the discrete time
domain where ${\glD}^\alpha f$ is approximated 
by its discrete time variant
\begin{align}
 ({\glDD}^\alpha x)_k =  
   \tfrac{1}{h^\alpha}
   \sum_{i=0}^{\infty}(-1)^i \tbinom{\alpha}{i}x_{k-i},
   \label{eq:GL_discrete}
\end{align}
which is in turn approximated by a discrete operator with 
finite memory $\nu$
\begin{align}
 ({\glDD}^\alpha_\nu x)_k =  
   \tfrac{1}{h^\alpha}
   \sum_{i=0}^{\nu}(-1)^i \tbinom{\alpha}{i}x_{k-i},
   \label{eq:GL_discrete_truncated}
\end{align}
which is proven to have bounded error with respect to 
$({\glDD}^\alpha f)_k$~\citep{SopSar16}.

%
%
\textit{Numerical integration methods:}
Fractional-order initial value problems can be solved
with various numerical methods such as the Adams-Bashforth-Moulton 
predictor-corrector (ABMPC) method~\citep{Zayernouri20161} and 
fractional linear multi-step methods (FLMMs)~\citep{Garrappa201596}.
These methods are only suitable for a system of FDE's in the form
\begin{subequations}\label{eq:FDE_basic_form}
\begin{align}
 {_{c}\D}^\gamma x(t) &= f(t,x(t)),\\
 x^{(k)}(0) &= x_{0,k},\, k = 1,\dots,m-1
\end{align}
\end{subequations}
where $\gamma$ is a rational, and $m = \lceil \gamma \rceil$.

In order to bring~\eqref{eq:pkModelAmiodarone} in this form, 
we need to find a rational approximation of two derivatives, $1-\alpha$ and $1$. 
If we can find a satisfying rational approximation of $1 - a\approxeq p/q$,
then first order derivative follows trivially.
Now,~\eqref{eq:pkModelAmiodarone} can be written as 
\begin{subequations}\label{eq:pkModelAmiodarone_system}
\begin{align}
    _{c}\D^{\gamma} x_0 &= x_1 \\ 
    _{c}\D^{\gamma} x_1 &= x_2 \\
    &\ \,\vdots \nonumber \\
    _{c}\D^{\gamma} x_{q-1} &= -(k_{12}+k_{10})x_0 {+} k_{21} x_{q+p}(t) {+} u\\
    _{c}\D^{\gamma} x_{q} &=  x_{q+1}(t)  \\
    &\ \, \vdots \nonumber \\
    _{c}\D^{\gamma} x_{2q-1} &=  k_{12}x_0 - k_{21}x_{q+p}
\end{align}
\end{subequations}
subject to $x_0 = A_1(0)$, $x_{q} = A_2(0)$ and $x_i(0) = 0$ for $i\neq 1, q$,
and $\gamma=1/q$.
This system is in fact a linear fractional-order 
system for which analytical solutions are available~\citep{Kaczorek2011}.

The number of states of system~\eqref{eq:pkModelAmiodarone_system} is $2q$, 
therefore, the rational approximation should 
aim at a small $q$. In our case $1-\alpha= 0.413$ can be written 
as $413/1000$, but then we would need to simulate a fractional-order
system with $2000$ states. Instead $1-\alpha$ can be approximated by 
$19/46$ with error $0.413-19/46 = -4.3478\cdot 10^{-5}$ or 
$216/523$ with error $0.413-216/523 = -1.912\cdot 10^{-6}$.
Such approximations can be obtained by means of continued fractions 
expansions of $1-\alpha$.
Yet another reason to choose small $q$ is that small values of 
$\gamma=1/q$ render the system hard to simulate numerically.

%
%
\textit{Adams-Bashforth-Moulton predictor-corrector (ABMPC):} 
Methods of the ABMPC type  have been generalized 
to solve fractional-order systems.
The basic concept is to evaluate $({\rlI}^\gamma f)(t,x(t))$ by 
approximating $f$ with appropriately selected polynomials.
Solutions of~\eqref{eq:FDE_basic_form} satisfy the following
integral representation
\begin{align}
 x(t) = \sum_{k = 0}^{m-1} x_{0,k}\frac{t^k}{k!} + ({\rlI}^\gamma f)(t,x(t)),
\end{align}
where the first term on right hand side will be denoted with $T_{m-1}(t).$
The integral on the right hand side of the previous equation can 
be approximated, using an uniformly spaced grid $t_n = nh$, by
$
  \frac{h^\gamma}{\gamma(\gamma + 1)}\sum_{j = 0}^{n+1} a_{j,n+1} f(t_j) 
$ for suitable coefficients $a_{j,n+1}$~\citep{Diethelm2002}.
The numerical approximation of the solution of~\eqref{eq:FDE_basic_form} 
is 
\begin{subequations}
\begin{align}\label{eq:corrector_equation}
 x(t_{n+1}) &= T_{m-1}(t_{n+1}) 
             + \tfrac{h^\gamma}
                     {\Gamma(\gamma+2)}f(t_{n+1},x_p(t_{n+1}))\notag\\
            &+ \sum_{j=1}^n a_{j,n+1}f(t_j,x(t_j)). 
 \end{align}
The equation above is usually referred to as the corrector formula and 
$x_p(t_{n+1})$ is given by the predictor formula
\begin{align}
 x_P(t_{n+1}) = T_{m-1}(t_n) &+ \tfrac{1}{\Gamma(\gamma)} 
                 \sum_{j = 0}^{n} b_{j,n+1}f(t_j,x(t_j).
\label{eq:predictor_equation}
\end{align}
\end{subequations} 
Unfortunately, the convergence error of ABMPC when $0 < \gamma < 1$
is $\mathcal{O}(h^{1+\gamma})$, therefore, a rather small step size 
$h$ is required to attain a reasonable aprpoximation error.
A modification of the basic predictor-corrector method with more 
favorable computational cost is provided in~\citep{Garrappa_pc_stability} 
for which the MATLAB implementation \texttt{fde12} is available.

\textit{Lubich's method:} Fractional linear multistep 
methods (FLMM)~\citep{Lubich} are a generalization of 
linear multistep methods (LMM) for ordinary differential 
equations. The key idea is to approximate the 
Riemann-Liouville fractional-order integral operator~\eqref{eq:RL_frac_integral}
with a discrete convolution, called \textit{convolution quadrature}, as
\begin{align}
 (\rlI_h^\gamma f)(t) 
   \approxeq h^\gamma\sum_{j=0}^{n} \omega_{n-j}f(t_j) 
   + h^\gamma \sum_{j=0}^{s}w_{n,j}f(t_j),
 \label{eq:flmm_conv_quadrature}
\end{align}
for $t_j = jh, h > 0$ where starting ($w_{n,j}$) and quadrature 
weights ($\omega_n$) are independent of $h$. 
%
Surprisingly, the latter weights can be constructed from any linear 
multistep method for arbitrary fractional order $\gamma$~\citep{Lubich}.
Furthermore, FLMM constructed this way will inherit the same 
convergence rate and at least the same stability properties as the 
original LMM method~\citep{Lubich_flmm_abel_voltera}.

Here we use the MATLAB implementation 
\texttt{flmm2}~\citep{Garrappa_software} which is based 
on~\citep{Garrappa_trapezoidal}.
However, the method does not perform well for small $\gamma$.
In our case study simulations we have found that values smaller than  
$0.1$ give poor results and often do not converge. However, when we
used a more crude approximation of the original system with $\gamma = 1/5$, 
\texttt{flmm2} method was outperforming \texttt{fde12} in terms of accuracy and 
has shown excellent stability properties with respect to bigger step size $h$.

\subsubsection{Numerical inverse Laplace.}
Several numerical inverse Laplace methods provide an approximation of
\begin{align}
  f(t) = \lim_{T \rightarrow \infty}\frac{1}{2\pi i} \int_{\sigma - iT}^{\sigma + iT} \e^{st}F(s)\d s.
  \label{eq:inv_laplace}
\end{align}
for a given transfer function $F(s)$.
Numerical methods can be used to directly evaluate the 
inversion integral~\eqref{eq:inv_laplace} for non-rational transfer functions. 
One of the most popular methods is to convert the 
inversion integral into a Fourier transform
and then approximate it by a Fourier series via trapezoid 
rule~\citep{Hoog_invlap}. These methods are considered to be precise for a broad class of function,
but computationally demanding. An implementation of the 
above method is freely available on line~\citep{Hollenbeck}.


A somewhat different approach is taken by \cite{valsa_invlap},
where authors approximate $\e^{st}$, the kernel of the 
inverse Laplace transformation, by $\tfrac{\e^{st}}{1 + e^{-2a}e^{2st}}$
and choose $a$ appropriately so as to achieve an accurate inversion.

In general, numerical inversion methods can achieve high precision, but 
they are not suitable for control design purposes, especially
for optimal control problems. Moreover, different methods are suited
for various types of problems.
An overview of the most popular inversion methods used in 
engineering practice is given in~\citep{invlap_comp}.

\subsection{Assessment}
In order to assess the accuracy of each approximation 
method presented above we introduce the following error
indices
\begin{align}
 \|e_i\| &= \sqrt{\int_0^{T_s} e_i(\tau)^2 \mathrm{d}\tau}, \\
 \|e_i\|_\infty &=  \max_{t\in [0,T_s]} e_i(t ), 
\end{align}
where $T_s$ is a fixed simulation time and $e_i(t)$ is the 
difference between the approximate response of the system
$\hat{A}_i$ and the one estimated by the inverse Laplace method
of~\cite{valsa_invlap} with $a=11$ which is considered to be the most accurate.

\begin{figure}[h]
\centering
\includegraphics[scale = 0.42]{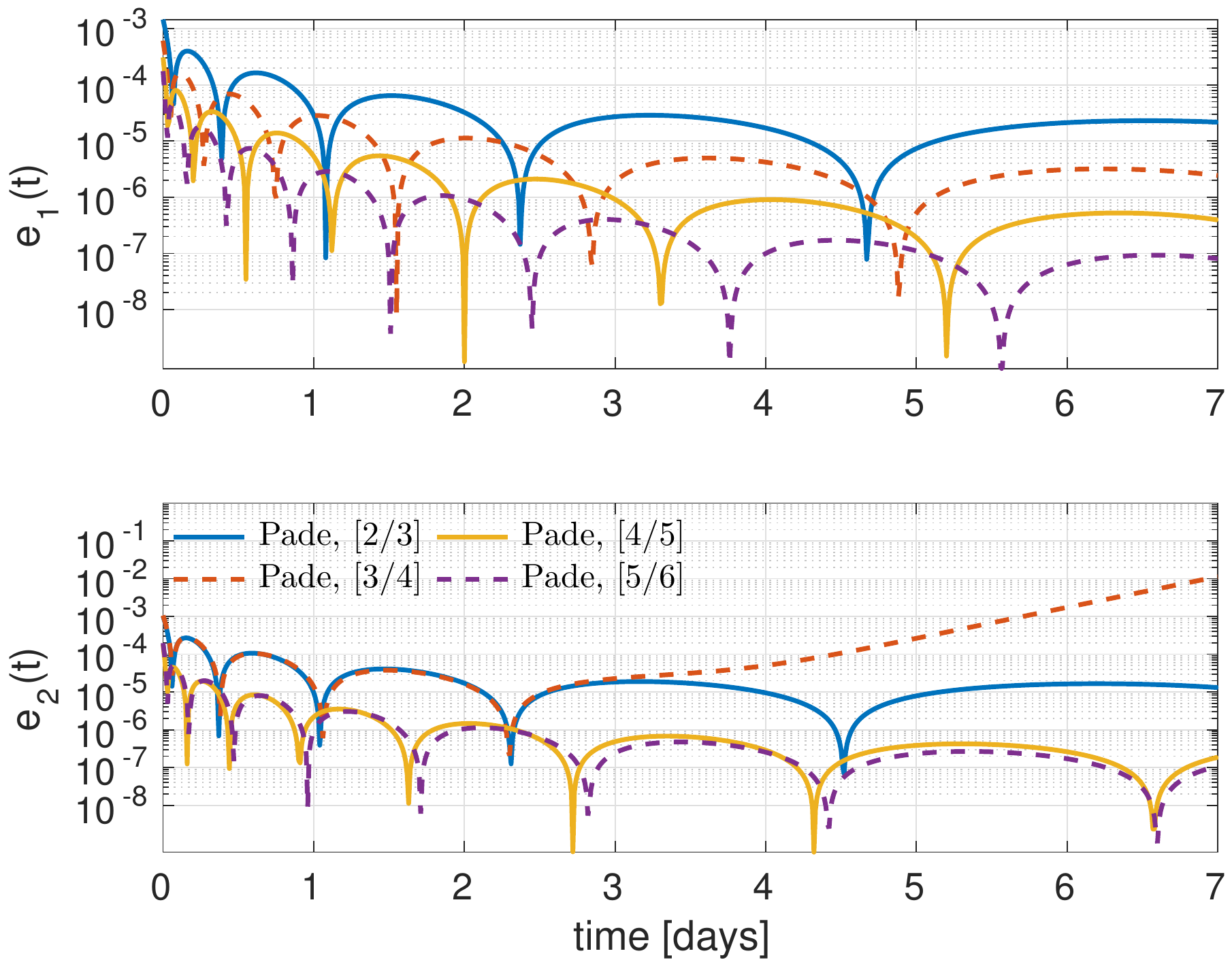}
\caption{Absolute error of Pad\'e  approximation 
         for various transfer function orders.}
\label{fig:pade}         
\end{figure}

\begin{figure}[h]
\centering
\includegraphics[scale = 0.42]{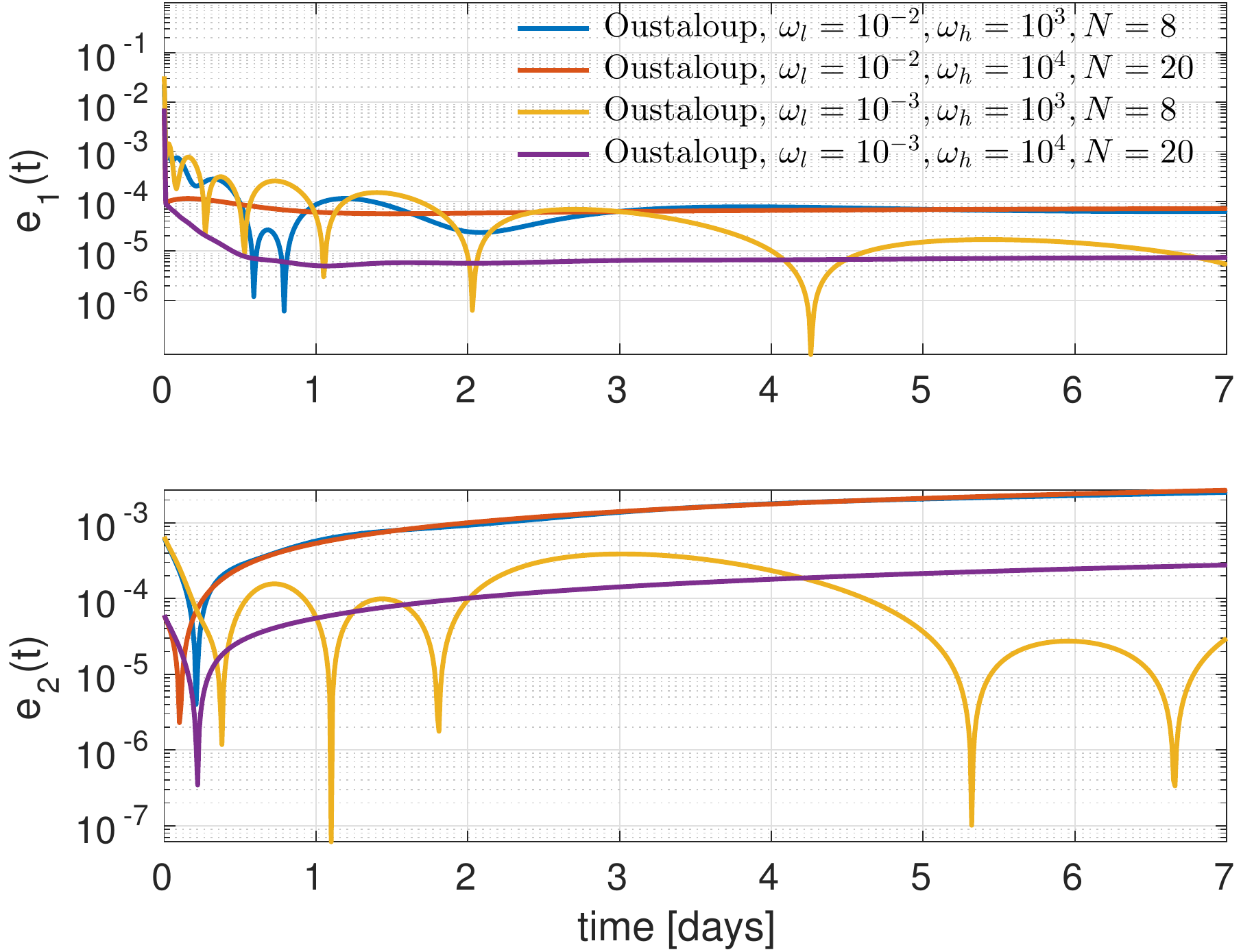}
\caption{Absolute error of Oustaloup approximation for various parameters.}
\label{fig:oustaloup}
\end{figure}

\begin{figure}[h]
\centering
\includegraphics[scale = 0.42]{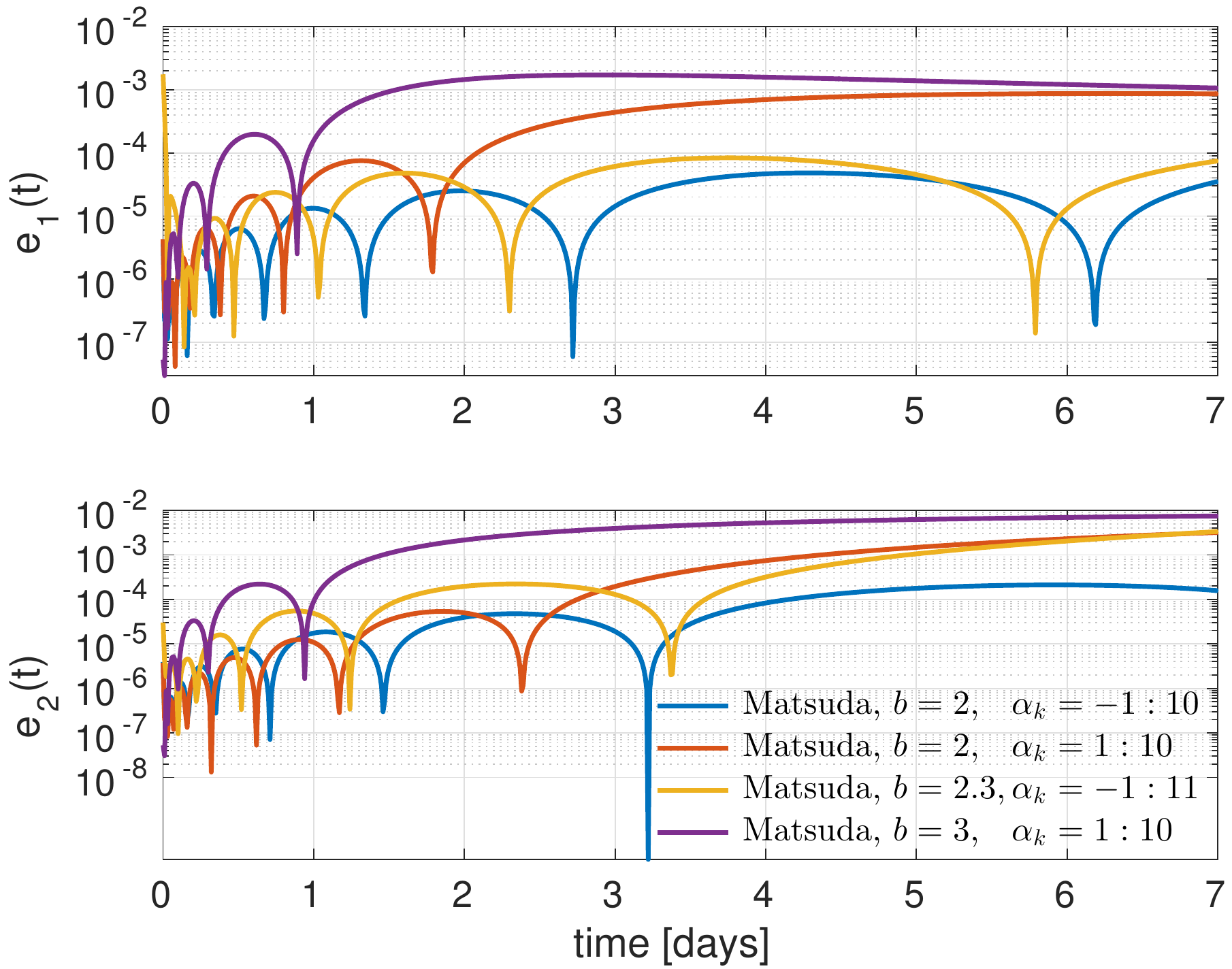}
\caption{Absolute error of Matsuda-Fujii approximation 
         for various parameters.}
\label{fig:matsuda-fujii}         
\end{figure}

The pharmacokinetic profile following a single \textit{i.v.} bolus
dose is shown in Figure~\ref{fig:open-loop-sim}.
In Figures~\ref{fig:pade},~\ref{fig:oustaloup} and~\ref{fig:matsuda-fujii}
we show the modeling errors $e_i(t)$ and in Tables~\ref{table:pade},%
~\ref{table:oustaloup} and~\ref{table:matsuda-fujii} we show the 
corresponding total errors. It seems that adequately high precision
can be achieved with these methods. However, approximations in the $s$-domain 
are not suitable for constrained systems since there is no theoretical
bound on the approximation error in the time domain.

\begin{figure}[h]
\centering
\includegraphics[scale = 0.42]{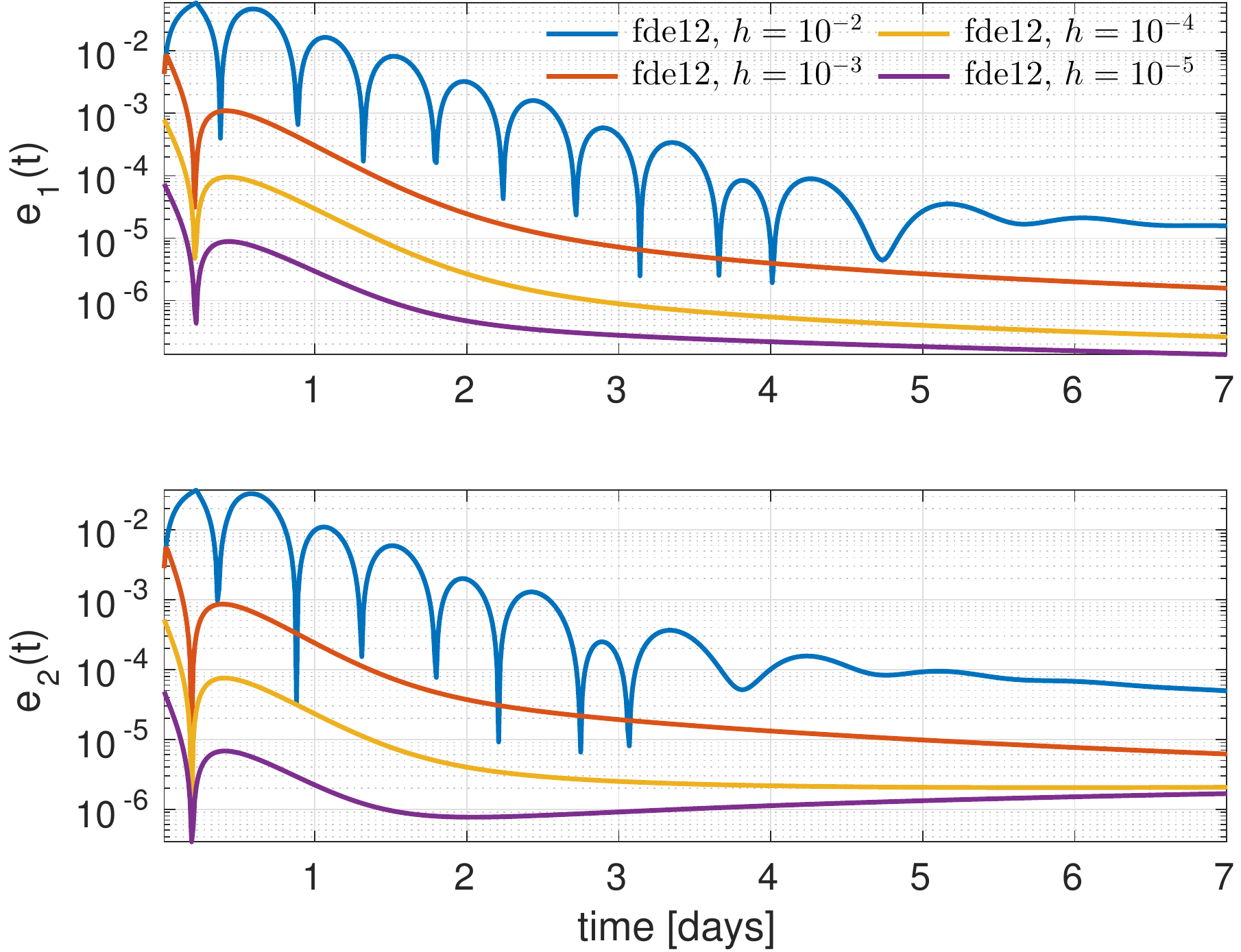}
\caption{Absolute error of fractional differential equation solution 
         by means of fde12 method where fractional order $1 - a = 0.413$ is 
         approximated by $19/46$. }
\label{fig:error-fde}        
\end{figure}

\begin{figure}[h]
\centering
\includegraphics[scale = 0.42]{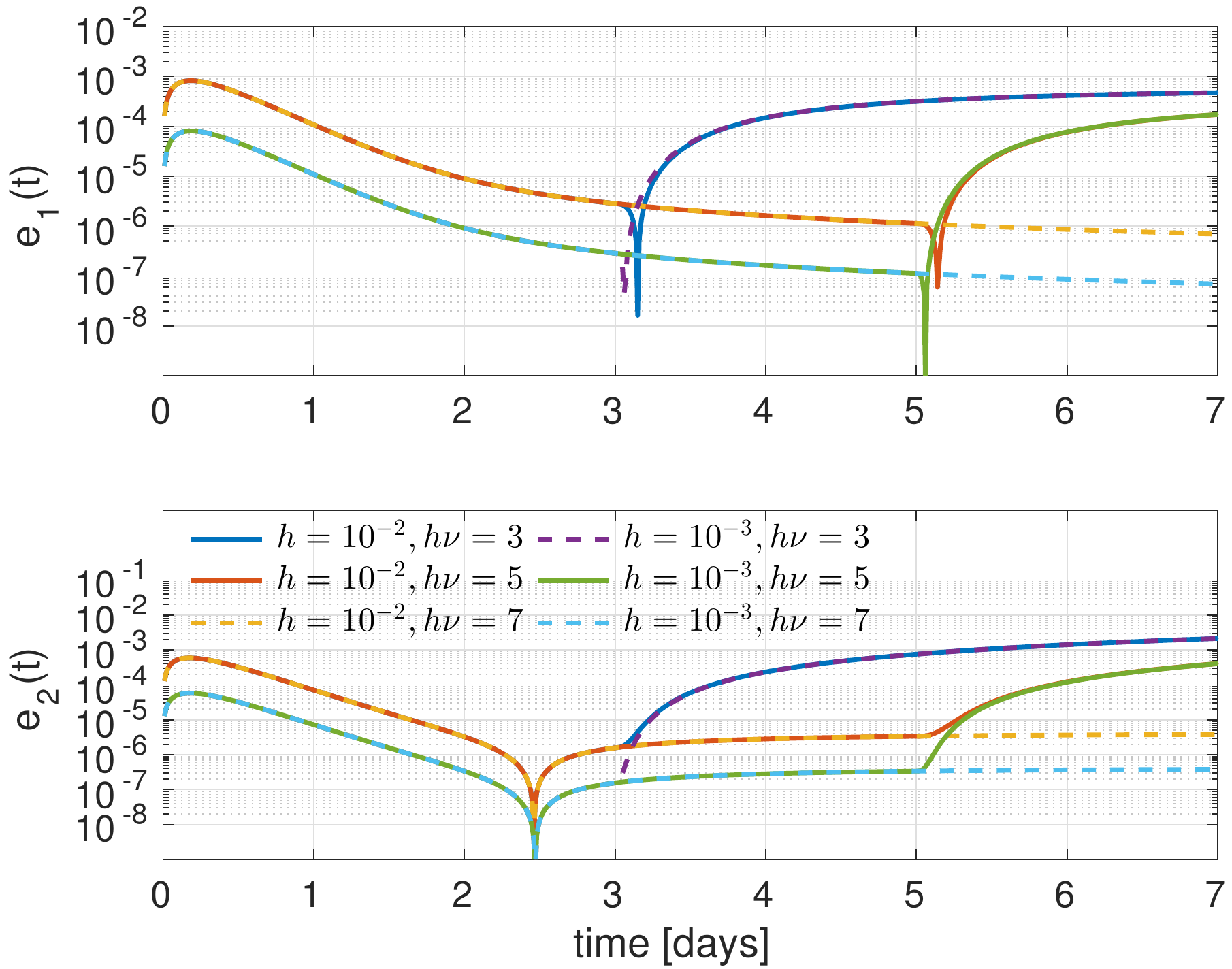}
\caption{Absolute error of GL method for various 
         step sizes and history lengths.}
\label{fig:error-gl}
\end{figure}

The errors of \texttt{fde12} are presented in Figure~\ref{fig:error-fde}
and Table~\ref{table:error-fde}. Using a step size as small as $h=10^{-5}$
\texttt{fde12} achieves an approximation error which is uniformly
lower than $10^{-4}$. In Figure~\ref{fig:error-gl} and Table~\ref{table:error-gl}
we show the approximation errors of the Gr\"unwald-Letnikov 
method. It can be seen that the use of a long history 
is more important for the attainment of high precision compared to
a small step size.
Most likely as a result of the low value of $\gamma=0.0217$
in~\eqref{eq:pkModelAmiodarone_system}, \texttt{flmm2}
failed to produce reasonable approximations. We were 
in fact only able to produce moderately accurate 
approximations using $1-\alpha\approxeq 2/5$ where $\gamma=0.2$.

%
%

%
%

%
%

%
%
%
\setlength{\tabcolsep}{0.4em} 
\begin{table}[h]
\caption{Errors using the Pad\'e approximation.}
\label{table:pade}
\centering
\begin{tabular}{c c c c c c c}
 Order & $\|e_1\|$  & $\|e_2\|$ & $\|e_1\|_\infty$  & $\|e_2\|_\infty$\\ \hline  
   $[2/3]$ 	&$2.833\cdot 10^{-4} $ &$1.907\cdot 10^{-4} $& $0.0015$ & $0.001$\\
   $[3/4]$ 	&$1.105\cdot 10^{-4}$ &$0.0059$& $6.094\cdot 10^{-4}$ & $0.0113$ \\
   $[4/5]$ 	&$4.514 \cdot 10^{-5}$ &$2.406 \cdot 10^{-5}$& $3.076\cdot 10^{-4}$ & $1.774\cdot 10^{-4}$ \\
   $[5/6]$ 	&$2.327 \cdot 10^{-5}$ &$2.685 \cdot 10^{-5}$& $1.752\cdot 10^{-4}$ & $1.976\cdot 10^{-4}$\\
\end{tabular}
\end{table}
 
\begin{table}[h]
\caption{Errors using the Oustaloup approximation.}
\label{table:oustaloup}
\centering
\begin{tabular}{c c c c c c c}
 $\omega_b$ & $\omega_h$& $N$ & $\|e_1\|$  & $\|e_2\|$ & $\|e_1\|_\infty$  & $\|e_2\|_\infty$\\ \hline  
 $10^{-2}$  & $10^3$ & $8$  & $0.0023$   & $0.0043$& $0.0228$ & $0.0025$ \\
 $10^{-2}$  & $10^4$ & $20$ & $5.744\cdot 10^{-4}$   & $0.0045$& $0.0054$ & $0.0027$\\
 $10^{-3}$  & $10^3$ & $8$  & $0.0033$   & $5.084\cdot 10^{-4}$& $0.0332$ & $6.555\cdot 10^{-4}$\\
 $10^{-3}$  & $10^4$ & $20$ & $7.451\cdot 10^{-4}$   & $4.597\cdot 10^{-4}$& $0.0074$ & $2.765\cdot 10^{-4}$\\ 
\end{tabular}
\end{table}

\begin{table}[h]
\caption{Errors using the Matsuda-Fujii 
         approximation with $s_k = \beta^{\alpha_k}$.}
\label{table:matsuda-fujii}         
\centering
\begin{tabular}{c c c c c c c}
 $\beta$ & $\alpha_k$& $\|e_1\|$  & $\|e_2\|$ & $\|e_1\|_\infty$  & $\|e_2\|_\infty$ \\ \hline  
   $2$ 	 &$-1:10$ &$7.01\cdot 10^{-5}$  &$3.162\cdot 10^{-4}$& $4.82 \cdot 10^{-5}$ & $2.111 \cdot 10^{-4}$ \\
   $2$ 	 &$1:10$  &$0.0016$  &$0.0036$& $0.0009$ & $0.0032$\\
   $2.3$ &$-1:11$ &$0.0002$  &$0.0032$& $0.0018$ & $0.003$\\
   $3$	 &$1:10$  &$0.0034$  &$0.0127$& $0.0017$ & $0.0075$\\
\end{tabular}
\end{table}


%
%

%
%

\begin{table}[h]
\caption{Approximation errors of \texttt{fde12}}
\label{table:error-fde}
\begin{center}
\begin{tabular}{c c c c c c c}
  $h$ & $\|e_1\|$  & $\|e_2\|$ & $\|e_1\|_\infty$  & $\|e_2\|_\infty$\\ \hline  
 $10^{-2}$  &$0.0333$ &$0.0223$ &    $0.0584$  &  $0.0368$ \\ 
 $10^{-3}$  &$0.0020$ &$0.0013$ &   $0.0083$  &  $0.0055$\\ 
 $10^{-4}$  &$1.782\cdot 10^{-4}$ & $1.157\cdot 10^{-4}$  & $7.982\cdot 10^{-4}$ & $5.204\cdot 10^{-4}$ \\ 
 $10^{-5}$  &$1.669 \cdot 10^{-5}$  &$1.113\cdot 10^{-5}$ & $7.412\cdot 10^{-4}$ & $4.824\cdot 10^{-4}$\\
\end{tabular}
\end{center}
\label{table:error-gl}
\end{table}

\begin{table}[h]
\caption{Approximation errors of the Gr\"unwald-Letnikov method}
\begin{center}
\bgroup
\def\arraystretch{1}
\begin{tabular}{c  c  c  c  c  c  } 
 $h$         & $h\nu$   & $\|e_1\|$              & $\|e_2\|$               &  $\|e_1\|_\infty$      &  $\|e_2\|_\infty$\\ \hline 
 $10^{-2}$   & $3$      &$8.165\cdot10^{-4}$     &$0.0022$                 & $8.223\cdot10^{-4}$    &   $0.0021$ \\
 $10^{-2}$   & $5$      &$5.302\cdot10^{-4} $    &$4.571\cdot 10^{-4} $    & $8.223 \cdot 10^{-4}$  &   $5.911\cdot 10^{-4}$ \\
 $10^{-2}$   & $7$      &$5.126\cdot 10^{-4} $   &$3.640\cdot 10^{-4} $    & $8.223 \cdot 10^{-4}$  &   $5.911\cdot 10^{-4}$ \\
 $10^{-3}$   & $3$      &$6.390\cdot 10^{-4}$     &$0.0022$                 & $4.735\cdot 10^{-4}$   &   $0.0021$    \\
 $10^{-3}$   & $5$      &$1.454\cdot 10^{-4}$    &$2.760 \cdot 10^{-4} $   & $1.727\cdot 10^{-4}$   &   $4.057\cdot 10^{-4}$ \\
 $10^{-3}$   & $7$      &$5.060 \cdot 10^{-5} $  &$3.594 \cdot 10^{-5}$    & $8.067\cdot 10^{-4}$   &   $5.797\cdot 10^{-4}$ \\ 
\end{tabular}
\egroup
\end{center}
\end{table}

\section{Administration scheduling}
In this section we address the problem of administration scheduling.
Our objective is to devise an administration schedule --- a sequence 
of dosages --- so that the concentration of Amiodarone in the tissues
is close to a desired value, while the concentration in both compartments
never exceeds certain safety limits.
We also have to account for a limit on allowed drug dose at each time instant. 
In doing so, we must assume that we are not able to measure drug concentrations 
during the treatment. All of these requirements and constraints 
can be elegantly integrated within the framework of constrained optimal control.

\subsection{Optimal control}
In this section we describe the optimal control problem formulation.
We start by discretizing~\eqref{eq:pkModelAmiodarone} with a sampling 
time $t_c$ yielding
\begin{align}
t_c^{-1}(x_{k+1}-x_{k}) = Ax_k + F \glDD^{1 - a}_{\nu} x_k + Bu_k 
\label{eq:asdf-qwerty} 
\end{align}
where $x_k  = \left[ A_1( kt_c) \,\, A_2( kt_c) \right]'$.
The left hand side of~\eqref{eq:asdf-qwerty} corresponds to the 
forward Euler approximation of the first-order derivative, and we shall refer to 
$t_c=\unit[10^{-2}]{days}$ as the \textit{control sampling time}. 
Matrices $A,F$ and $B$ are
\begin{align}
A = \begin{bmatrix}
    -(k_{12} + k_{10}) & 0 \\
      k_{21}        & 0 \\
\end{bmatrix},
F = \begin{bmatrix}
     0 & k_{21}  \\
     0 &-k_{21}  \\
\end{bmatrix},
B = \begin{bmatrix}
     1\\
     0\\
\end{bmatrix}.
\end{align}
The discrete-time dynamic equations of the system can now be stated as 
\begin{align}
    x_{k+1} = x_k 
            + t_c \bigg(
               Ax_k + \frac{F}{t_c^{1-a}} \sum_{j = 0}^{\nu} c_j^{1-a}x_{k-j} 
                    + Bu_k
                 \bigg),
    \label{eq:dynamic_equation}
\end{align}
where $c_j^\alpha = (-1)^j\binom{\alpha}{j} $.
By augmenting the system with past values as 
$\tilde{x}_k = (x_k,x_{k-1},\dots,x_{k-\nu + 1})$  we can 
rewrite~\eqref{eq:dynamic_equation} 
as a finite-dimension linear system
\begin{align}
  \tilde{x}_{k+1} = \hat{A} \tilde{x}_k + \hat{B}u_k.
  \label{eq:dynamic_equation_LTI} 
\end{align}
Matrices $\hat{A}$ and $\hat{B}$ are straightforward to derive and are given in~\citep{SopSar16}.
The therapeutic session will last for $N_d = N t_c = \unit[7]{days}$ 
in total, where $N$ is called the \textit{prediction horizon}.
It is not realistic to administer the drug to  the patient too frequently, 
so we assume that the patient is to receive their treatment every $t_d=\unit[0.5]{days}$.
The administration schedule must ensure that the concentration of 
drug in all compartments never exceeds the minimum toxic concentration 
limits while tracking the prescribed reference value as close as possible.
To this aim we postulate the 
following constrained optimal control problem.
\begin{subequations}\label{eq:osocp}
 \begin{align}
   \min_{ \{u_0,\ldots, u_{N_d-1}\}}J &= \sum_{k=0}^{N_d/t_c + 1}  (x_{\mathrm{ref},k} - x_k)' Q (x_{\mathrm{ref},k} - x_k) \\
 \intertext{subject to}
  \tilde{x}_{k+1} &= \hat{A} \tilde{x}_k + \hat{B} u_j, \text{ for } \ k t_c = j t_d\\
  \tilde{x}_{k+1} &= \hat{A} \tilde{x}_k, \text{ otherwise}\\
  0 &\le x_{k} \le 0.5 \\
  0 &\le u_{j} \le 0.5 \\
  \intertext{for $k=0,\ldots, N$; $j=0,\ldots, N_d-1$.} \nonumber
  \end{align}
\end{subequations} 
In the above formulation $x_{\mathrm{ref},k}$ is the desired drug concentration at time $k$ 
and operator $'$ denotes vector transposition.
Any deviation from set point is penalized by weight matrix $Q = \mathrm{diag}([0\;1])$. 
Note that we are tracking only the second state. 
Our underlying GL model has a relative history of $t_c \nu =\unit[5]{days}$.
Optimal drug concentrations are denoted by $u^\star_k$, for $k = 0,\dots, N_d-1$ and they correspond
to dosages administered intravenously at times $kt_d$. In the optimal control formulation we have implicitly
assumed that $t_d$ is an integer multiple of $t_c$, which is not restrictive since $t_c$ can be chosen arbitrarily.
Finally, we can recognize that problem~\eqref{eq:osocp} is a standard quadratic problem that can be readily solved. 
In our simulations, we have used YALMIP~\citep{YALMIP} to model the problem and MOSEK~\citep{mosek} 
as the underlying solver to calculate the solution.

\subsection{Simulations}
To argue for the soundness and the applicability of our approach in real-world scenarios, 
we will apply the optimal drug dosage schedule to a more precise model than~\eqref{eq:dynamic_equation_LTI}.
For this purpose we will use \texttt{fde12} solver. As evident from the results in the previous section,
for sufficiently small solver time $h_{sol}$, we can have a realistic simulation of the system. 

\begin{figure}[h]
\centering
\includegraphics[scale = 0.42]{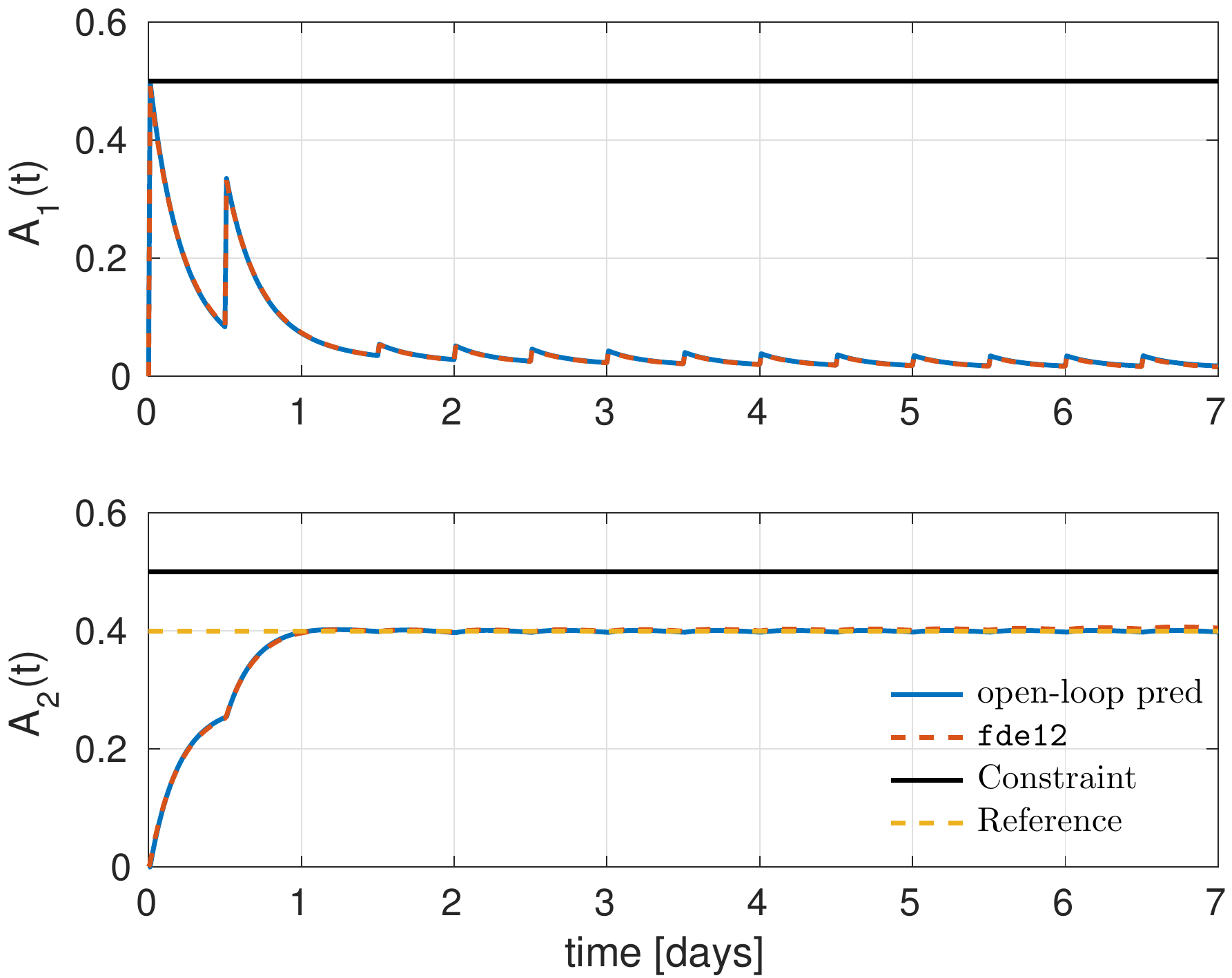}
\caption{Open loop control of drug administration with a fixed scheduling. 
Step size for \texttt{fde12} was set to $h = 10^{-5}$.}
\label{fig:sim_one_patient_response}
\end{figure}
\begin{figure}[h]
\centering
\includegraphics[scale = 0.42]{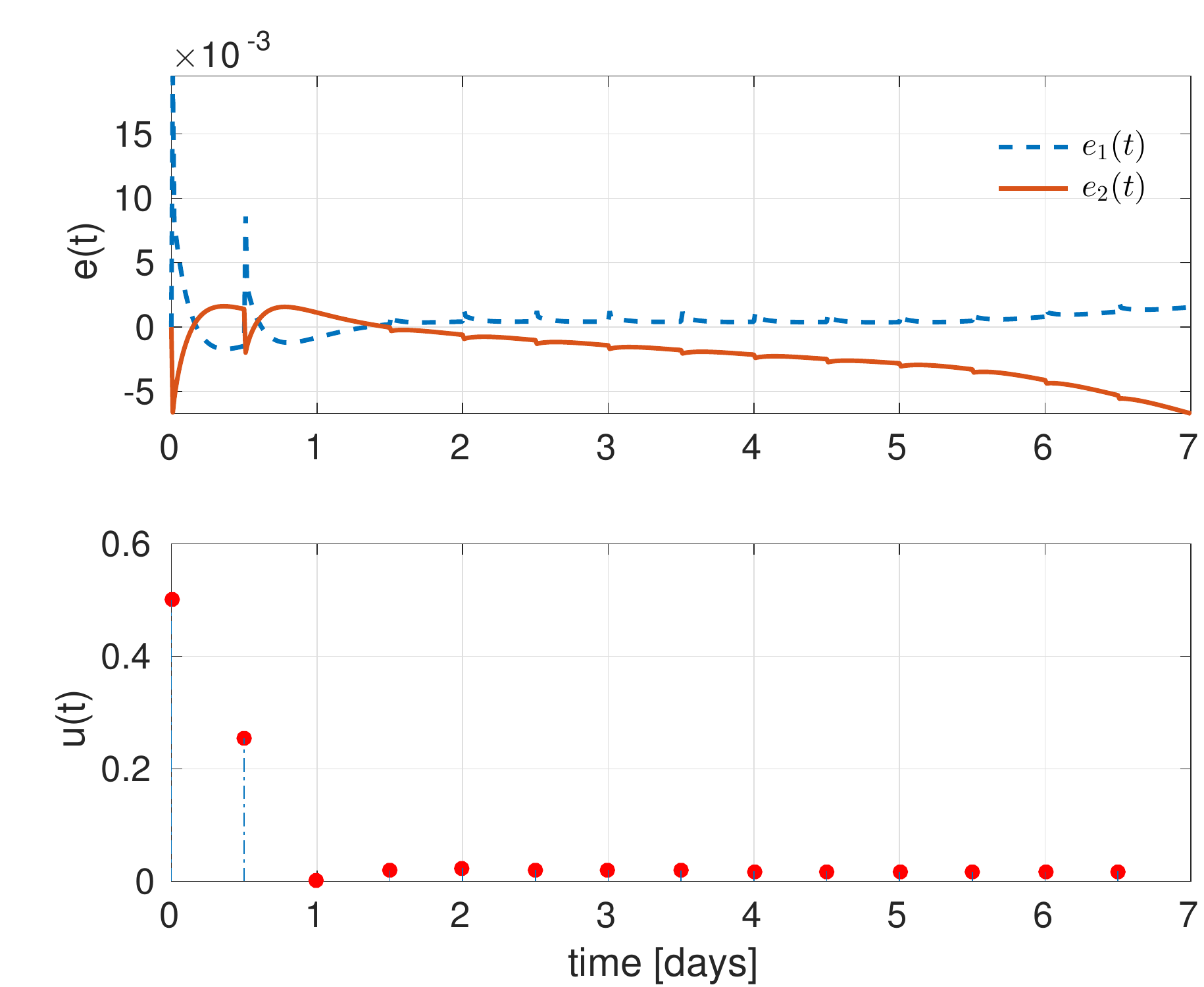}
\caption{(Up) errors of predicted states via GL model against \texttt{fde12}. 
(Down) optimal administration sequence.}
\label{fig:sim_one_patient_error}
\end{figure}

After solving the optimal control problem we applied the optimal sequence to the FDE simulator \texttt{fde12}.
Results are shown in Figures~\ref{fig:sim_one_patient_response} and~\ref{fig:sim_one_patient_error}. 
It can be seen that open loop predictions of the GL model and \texttt{fde12} simulation show high agreement.
This should not be surprising considering that all of the parameters describing the patient are nominal
and, additionally, shows that the GL model is of good quality.

\begin{figure}[h]
\centering
\includegraphics[scale = 0.42]{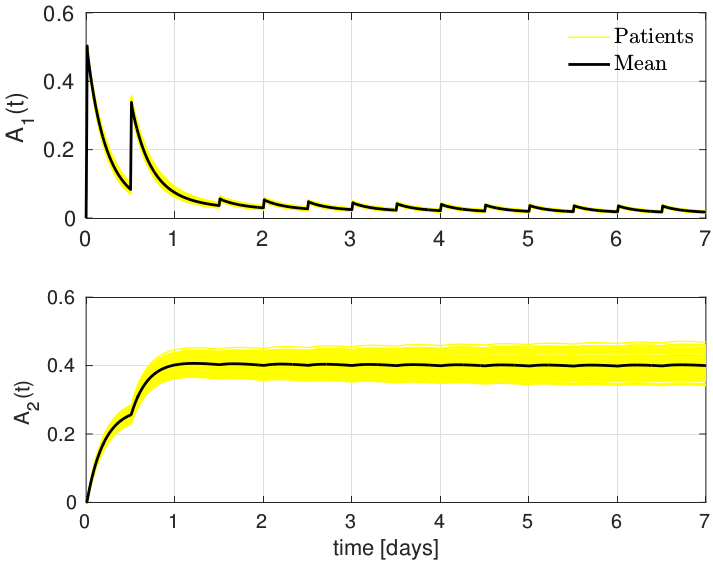}
\caption{Fixed schedule drug administration for a population of $100$ patients. 
Step size for \texttt{fde12} was set to $h = 10^{-5}$.}
\label{fig:sim_multiple_response}
\end{figure}

Next, we consider multiple patients that are characterized by perturbed parameters.
We will simulate the behaviour of $100$ different patients with multiplicative perturbations on parameters
$k_{10},k_{12}$ and $k_{21}$. Each parameter is multiplied by a constant 
drawn from a uniform distribution on the interval $\left[0.85, 1.15\right]$. Moreover, 
the order of the fractional derivative $\alpha = 1 - \hat{p}/q$ is given by a random choice of $\hat{p}$ from 
a set of discrete values $\{17,18,20,21\}$, each one having the same probability. Denominator $q$ is
fixed at $q = 46$, while the nominal numerator is $p = 19$. 
We simulate each patient by applying the same optimal drug scheduling sequence that was computed for the nominal one, that is, 
without online information about the state of each patient.
Results are shown in Figure~\ref{fig:sim_multiple_response}.

\section{Conclusions}
This paper gives an overview of the state of the art in numerical methods for the 
simulation of fractional-order systems along with validation 
results on a fractional-order pharmacokinetic model taken from the literature. 
Results are shown regarding the solution of an open-loop optimal control problem 
for the administration of Amiodarone to a patient whose pharmacokinetic
parameters are assumed to be perfectly known. 
Additionally, we presented optimal control results for different patients whose pharmacokinetic parameters 
are not know perfectly. This way, we estimate and demonstrate how sensitive 
the administration scheduling is.

\bibliography{ifacconf}

\end{document}